\documentclass[oneside,notitlepage,12pt]{article}
\newif\iffarben
\farbentrue
\usepackage[utf8]{inputenc}
\usepackage[T2A,T1]{fontenc}
\usepackage{geometry}
\usepackage[colorlinks]{hyperref}
\usepackage{amsmath}
\usepackage{amssymb}
\usepackage[english,russian]{babel}

\usepackage{moreverb}
\usepackage{placeins}

\newcommand\BibTeX{{\rmfamily B\kern-.05em \textsc{i\kern-.025em b}\kern-.08em
T\kern-.1667em\lower.7ex\hbox{E}\kern-.125emX}}

\usepackage{amsfonts}
\usepackage{amsopn}
\usepackage{amstext}
\usepackage{amsthm}
\usepackage{enumitem}

\usepackage{calrsfs}
\usepackage{clock} 
\usepackage{ifsym}

\usepackage{tikz}
\usetikzlibrary{intersections}
\usetikzlibrary{calc}
\usetikzlibrary {decorations.pathmorphing, decorations.pathreplacing, decorations.shapes}
\usetikzlibrary {shapes.geometric}

\definecolor{L1u}{rgb}{0,0,0}
\definecolor{L2u}{rgb}{0,0.1,0}
\definecolor{L3u}{rgb}{0,0.2,0}
\definecolor{L4u}{rgb}{0,0.3,0}
\definecolor{L5u}{rgb}{0,0.4,0}
\definecolor{L6u}{rgb}{0,0.5,0}
\definecolor{L7u}{rgb}{0,0.6,0}
\definecolor{L8u}{rgb}{0,0.7,0}
\definecolor{L9u}{rgb}{0,0.8,0}
\definecolor{L10u}{rgb}{0,0.9,0}

\definecolor{L1d}{rgb}{0,0.05,0}
\definecolor{L2d}{rgb}{0,0.15,0}
\definecolor{L3d}{rgb}{0,0.25,0}
\definecolor{L4d}{rgb}{0,0.35,0}
\definecolor{L5d}{rgb}{0,0.45,0}
\definecolor{L6d}{rgb}{0,0.55,0}
\definecolor{L7d}{rgb}{0,0.65,0}
\definecolor{L8d}{rgb}{0,0.75,0}
\definecolor{L9d}{rgb}{0,0.85,0}

\definecolor{R1u}{rgb}{0,0,0}
\definecolor{R2u}{rgb}{0.1,0,0}
\definecolor{R3u}{rgb}{0.2,0,0}
\definecolor{R4u}{rgb}{0.3,0,0}
\definecolor{R5u}{rgb}{0.4,0,0}
\definecolor{R6u}{rgb}{0.5,0,0}
\definecolor{R7u}{rgb}{0.6,0,0}
\definecolor{R8u}{rgb}{0.7,0,0}
\definecolor{R9u}{rgb}{0.8,0,0}

\definecolor{R1d}{rgb}{0.05,0,0}
\definecolor{R2d}{rgb}{0.15,0,0}
\definecolor{R3d}{rgb}{0.25,0,0}
\definecolor{R4d}{rgb}{0.35,0,0}
\definecolor{R5d}{rgb}{0.45,0,0}
\definecolor{R6d}{rgb}{0.55,0,0}
\definecolor{R7d}{rgb}{0.65,0,0}
\definecolor{R8d}{rgb}{0.75,0,0}
\definecolor{R9d}{rgb}{0.85,0,0}

\definecolor{M1u}{rgb}{0.2,0.0,0.9}
\definecolor{M1d}{rgb}{0.2,0.0,0.9}

\usepackage{subcaption}
\providecommand{\cal}{\mathcal}
\renewcommand{\Bbb}{\mathbb}

\newenvironment{pf}{\begin{proof}}{\end{proof}}

\newcommand{\El}{{\cal{L}}}
\newcommand{\Err}{{\Bbb{R}}}

\newcommand{\lam}{{\lambda}}
\newcommand{\al}{\alpha}

\newcommand{\rest}{\restriction}
\newcommand{\loe}{\leqslant}
\newcommand{\goe}{\geqslant}
\newcommand{\subs}{\subseteq}
\newcommand{\nnempty}{\ne\emptyset}

\newcommand{\Int}{\operatorname{int}}

\newcommand{\conv}{\operatorname{conv}}

\newcommand{\oraz}{\qquad\text{and}\qquad}

\newcommand{\setof}[2]{\{#1\colon #2\}}
\newcommand{\sett}[2]{\{#1\}_{#2}}
\newcommand{\sn}[1]{\{#1\}} 
\newcommand{\dn}[2]{\{#1,#2\}} 
\newcommand{\map}[3]{#1\colon #2 \to #3} 
\newcommand{\img}[2]{#1 #2} 
\newcommand{\inv}[2]{{#1}^{-1}{#2}} 

\newcommand{\bX}{{\mathbb{X}}}

\newif\ifnotes
\notestrue
\ifnotes
\usepackage[capitalize,nameinlink]{cleveref}
\crefname{empty}{}{}
\crefname{equation}{}{}
\Crefname{figure}{Figure}{Figures}
\crefname{page}{page}{pages}
\Crefname{enumi}{}{}
\Crefname{subsection}{Subsection}{Subsections}
\def\theoremname{Theorem}%
\def\examplename{Example}%
\newtheorem{thm}{\theoremname}[section]
\newtheorem{tw}[thm]{Theorem}
\newtheorem{wn}[thm]{Corollary}
\newtheorem{lm}[thm]{Lemma}
\newtheorem{prop}[thm]{Proposition}
\theoremstyle{definition}
\newtheorem{df}[thm]{Definition}
\newtheorem{ex}[thm]{Example}
\newtheorem{pyt}[thm]{Question}
\theoremstyle{remark}
\newtheorem{uwgi}{Remark}

\theoremstyle{plain}
\theoremstyle{definition}

\else
\newtheorem{tw}{Theorem}[section]
\newtheorem{wn}[tw]{Corollary}
\newtheorem{lm}[tw]{Lemma}
\newtheorem{prop}[tw]{Proposition}

\theoremstyle{definition}

\newtheorem{ex}[tw]{Example}

\theoremstyle{remark}

\fi

\newcommand{\btw}{\flat}
\newcommand*\samethanks[1][\value{footnote}]{\footnotemark[#1]}
\newcommand{\ignoreSpellCheck}[1]{#1}

\title{Monotone mappings and lines}
\author{{\sc Wiesław Kubi\'s}\footnote{Research supported by grant 20-22230L (Czech Science Foundation).}\\ \\
{\small Institute of Mathematics,} {\small Czech Academy of Sciences}\\
\and
{\sc Janusz Morawiec}\thanks{Research supported by the University of Silesia  Mathematics Department (Iterative Functional Equations and Real Analysis program)}\\
{\small Institute of Mathematics, University of Silesia in Katowice, Poland}\\
\and
{\sc Thomas Zürcher}\samethanks\\
{\small Institute of Mathematics, University of Silesia in Katowice, Poland}
}
\date{\today\ \clocktime}

\begin{document}
\selectlanguage{english}
\maketitle

\abstract{We study betweenness preserving mappings (we call them \emph{monotone}) defined on subsets of the plane. Once the domain is a convex set, such a mapping is either the restriction of a homography, or its image is contained in the union of a line and a single point, or its image consists of five points, one of them being between two disjoint pairs of the other four points.
We also show that an open planar set cannot be mapped in a one-to-one monotone way into the real line. From this we deduce that a one-to-one monotone mapping from a convex planar set with nonempty interior is necessarily a partial homography. Finally, we prove that a set consisting of three pairwise non-parallel lines does not admit a one-to-one monotone mapping into the real line, while on the other hand a set consisting of three closed line segments intersecting at a single point does admit such a mapping.}

\noindent
{\bf \ignoreSpellCheck{MSC} (2020):} 52A10,
52C45, 
03E05. 


\maketitle
\tableofcontents
\section{Introduction}

In the abstract setting, a \emph{betweenness} is a ternary relation $\btw$ on a set $X$, satisfying some natural axioms, see~\cite{Vel93}, \cite{Sol84}, and~\cite{Kub02}.
There is quite a lot of literature on this topic, see e.g.\ the survey~\cite{Pamb}.
We are only interested in two types of betweennesses: \emph{Euclidean} and \emph{linear}.
Euclidean betweenness in a real vector space is the standard one: $x$ is between $a$ and $b$ (denoted by $\btw(a,x,b)$) if and only if $x = (1-\lam)a + \lam b$ for some $0\loe \lam \loe1$, namely, $x$ belongs to the line segment joining $a$ and $b$.
This generalizes to vector spaces over arbitrary ordered fields, however we are not going to explore this direction.
Linear betweenness is the one induced by a fixed linear ordering: $x$ is between $a$ and $b$ if and only if $a \loe x \loe b$ or $b \loe x \loe a$.
In the real line, the Euclidean and the linear betweenness are of course the same.

Given two betweenness structures $(X,\btw_X)$, $(Y,\btw_Y)$, a mapping $\map{f}{X}{Y}$ will be called \emph{monotone} if it preserves the betweenness, that is,
\begin{equation*}
\btw_X(a,x,b) \implies \btw_Y(f(a),f(x),f(b))
\end{equation*}
for every $a,x,b \in X$.
A common name is \emph{betweenness preserving} mapping (see e.g.~\cite{HM08}).
Note that if $\btw_X$~and~$\btw_Y$ are linear, then this notion coincides with the usual one: A mapping is monotone if and only if it is increasing or decreasing (in the non-strict sense). A monotone \emph{isomorphism} is a monotone bijection whose inverse is monotone.
Every monotone bijection between linearly ordered sets is an isomorphism, but this is not the case for Euclidean betweenness: Take $X = \{a,b,c\} \subs \Err^2$, where $a,b,c$ are not on a single line, and $Y = \{a',b',c'\}$, where $a',b',c'$ are collinear; then any bijection from $X$ to $Y$ is monotone, while its inverse is not monotone.
In fact, among subsets of the plane one can easily find those in which the Euclidean betweenness is \emph{discrete}, that is, $\btw(x,y,z)$ holds if and only if $y \in \dn x z$.
A typical example is the circle $S^1$. Another example is any set $S$ with the property that $S \cap L$ contains exactly two points for every line $L \subs \Err^2$. Such a set can be easily constructed by transfinite induction, knowing that each line has cardinality continuum and distinct lines intersect in at most one point. In particular, $S$ is isomorphic to $S^1$ and any bijection between those two sets is a monotone isomorphism.
This is a clear evidence that one should look at sets that are either convex or at least ``resemble'' convex sets, having sufficiently many collinear triples and perhaps satisfying some axioms involving intersections of line segments.

Affine transformations are obvious examples of monotone mappings between real vector spaces. Perhaps lesser known are projective transformations, defined on half-spaces only. We shall discuss them in Section~\ref{SectGoodDobri} below.

When it comes to arbitrary subsets of the real plane, as mentioned above, the induced betweenness may easily become discrete. Again, a good example is any circle. In particular, every mapping defined on a circle is automatically monotone.
When it comes to sets of small cardinality, one can concentrate on the plane, due to the following observation (a proof is given in the appendix).

\begin{prop}\label{propreducn2}
	Let $X$ be a subset of a real vector space $V$, and let $V_0$ be a $2$\nobreakdash-dimensional linear subspace of $V$. If the cardinality of $X$ is strictly less than the continuum, then there exists a linear map $\map L V V$ such that $\img L V = V_0$ and $L \rest X$ is a monotone isomorphism onto its image.
\end{prop}

We address the question of classifying monotone mappings on ``reasonable'' (e.g.\ convex or with nonempty interior) subsets of Euclidean spaces.
We shall concentrate on the plane, which is the smallest non-trivial case and many arguments concerning betweenness easily extend to higher dimensions.

Among the results, we obtain a classification of monotone mappings defined on convex planar sets (Theorem~\ref{THMgdyifgi} below). Namely, either the image is contained in a line plus a single point outside of it, or the image is a certain five-element set, or else the mapping is a partial homography. Additionally, if the domain is an open convex set then either the image is contained in a line or else the mapping is a partial homography. As a corollary we obtain that an injective monotone mapping defined on a convex planar set is necessarily a restriction of a homography. Recall that a homography $h$ (also called \emph{projective transformation}) is monotone on either of the two open half-planes whose union is the domain of $h$. Of course, if $h$ is affine then its domain is the whole plane.

In the literature one can find several works studying mappings which preserve collinearity, typically assuming injectivity and aiming to show that the mapping is affine or projective.
This is closely related to the fundamental theorem of affine geometry: Every bijection of a real vector space that maps lines onto lines is affine.
Removing the assumption of being one-to-one, one can easily find examples of pathological (discontinuous) monotone mappings, see Section~\ref{SectZlosciBad} below.
On the other hand, our Corollary~\ref{WNdnfosdogi} can be viewed as yet another variant of the fundamental theorem of affine (or rather projective) geometry, where ``preserving lines'' is replaced by ``preserving betweenness'', while at the same time allowing arbitrary convex sets instead of the full plane.

\section{The role of homographies}\label{SectGoodDobri}

What are the ``nice'' monotone isomorphisms between convex sets? The first thing that comes to our mind is: affine transformations. These are precisely those monotone mappings that preserve not only the betweenness but also the $1$\nobreakdash-dimensional barycentric coordinates of the point in an interval. Specifically, if $x = (1-\lam)a + \lam b$ with $\lam \in [0,1]$ and $f$ is affine, then $f(x) = (1-\lam)f(a) + \lam f(b)$.
No other monotone mappings have this property.
Namely, suppose $\map f A B$ preserves the $1$\nobreakdash-dimensional coordinates and assume $A$ is a nonempty convex set.
Then $f$ extends uniquely to an affine transformation $\map{\widetilde f}{V_A}{V_B}$, where $V_A$, $V_B$ denote the linear spans of $A$~and~$B$, respectively.

Recall that a \emph{homography} is an isomorphism of real projective spaces, that is, a bijection mapping lines onto lines.
Note that in the real projective space the betweenness becomes ambiguous, since lines are actually circles. Betweenness can only be defined locally. In any case, every affine isomorphism extends to the projective space, however there are other projective transformations.
Specifically, if $n>1$ is a natural number, then the general form of a projective transformation of $V = \Err^n$ is
\begin{equation*}
f(x) = \frac 1 {\alpha + (v | x)} T(x),
\end{equation*}
where $\al \in \Err$, $v \in V$, $T$ is an affine transformation, and $(v | x)$ denotes the usual (real) scalar product of $v$ and $x$.
We can actually treat the formula above as the definition of a projective transformation\footnote{Typically, \emph{projective transformation} is synonymous to \emph{homography}, but there seems to be no special name for the more general, not necessarily one-to-one, variant. That is why we have decided to call them projective transformations, while only the bijective ones will be called homographies.}.
Note that it is defined on the union of the two open half-spaces:
\begin{equation*}
H^+ = \setof{x \in V}{(v|x)+\alpha > 0} \oraz H^- = \setof{x \in V}{(v|x)+\alpha <0}.
\end{equation*}
In case $v=0$, either $H^+=V$ or $H^-=V$ and $f$ is affine.
It is easy to check that both $f \rest H^+$ and $f \rest H^-$ are monotone and $H^-$, $H^+$ are maximal convex sets on which $f$ is monotone.
A mapping $f$ defined on a subset of~$\Err^n$ will be called a \emph{partial projective transformation} (\emph{partial homography})
if there is a projective transformation (homography) extending $f$.

\begin{tw}\label{THMhomographies}
	Let $G \subs \Err^2$ be a convex set, let $\map f G {\Err^2}$ be a monotone mapping with the following property:
	\begin{enumerate}
		\item[] There exist $a,b,c, d \in G$ such that $f(d) \in \Int [f(a), f(b), f(c)]$.
	\end{enumerate}
	Then $f$ is the restriction of a unique homography.
	In particular, $f$ is a homeomorphism and, at the same time, a monotone isomorphism onto its image.
\end{tw}

\begin{pf}
	First, note that there exists a unique homography $g$ defined on a half-plane containing the triangle $\Delta := [a,b,c]$, satisfying $f \rest \{a,b,c,d\} = g \rest \{a,b,c,d\}$.
	This is a general fact in projective geometry, however a self-contained proof can be found as the one of~Lemma~2.14 in~\cite{ArtFloMil}.
	
	Now, using induction, we can see that $f$ and $g$ actually coincide on a dense subset~$D$ of~$\Delta$. Namely, we can choose $D$ as the minimal subset of $\Delta$ containing $a,b,c,d$, that is closed under intersections of line segments. More precisely, for every $x_0,x_1, y_0, y_1 \in D$, if $[x_0, x_1] \cap [y_0 , y_1]$ is a single point then this point belongs to $D$, and if there is a unique $z \in [x_0,x_1]$ such that $y_1 \in [y_0,z]$ then $z \in D$.
	The idea is explained in Figure~\ref{FiguraJedna}.
	\begin{figure}
    \begingroup
    \iffarben
    \else
    \colorlet{orange}{black}
    \colorlet{red}{black}
    \colorlet{cyan}{black}
    \colorlet{blue}{black}
     \fi
		\begin{tikzpicture}[scale=4]
			\coordinate (fa) at (0,0);
			\coordinate (fb) at (1,0);
			\coordinate (fc) at (.5,1);
			\draw (fa)--(fb)--(fc)--(fa);
			\coordinate (fq) at ($(fc)!0.35!(fb)$);
			\coordinate (u') at ($(fa)!0.35!(fb)$);
			\coordinate (v') at ($(fa)!0.75!(fb)$);
			\coordinate (u'') at ($(fa)!0.40!(fb)$);
			\coordinate (v'') at ($(fa)!0.70!(fb)$);
			\draw[ultra thick,orange] (u'')--(v'');
			\path[name path=fafq] (fa)--(fq);
			\path[name path=u'fc] (u')--(fc);
			\path[name path=v'fc] (v')--(fc);
			\path[name path=u'v'] (u')--(v');
			\coordinate (fd) at ($(fa)!0.70!(fq)$);
			\fill (fd) circle (0.34375pt);
			\node[below right] at (fd) {$d$};
			\fill (fa) circle (0.34375pt);
			\fill (fb) circle (0.34375pt);
			\fill (fc) circle (0.34375pt);
			\node[below] at (fa) {$a$};
			\node[below] at  (fb) {$b$};
			\node[above] at  (fc) {$c$};

			\begin{scope}[scale=1,xshift=-1.2cm]

				\draw[ultra thick,orange] (u'')--(v'');
			\end{scope}
		\end{tikzpicture}\hfill
		\begin{tikzpicture}[scale=4]
			\coordinate (fa) at (0,0);
			\coordinate (fb) at (1,0);
			\coordinate (fc) at (.5,1);
			\draw (fa)--(fb)--(fc)--(fa);
			\coordinate (fq) at ($(fc)!0.35!(fb)$);
			\coordinate (u') at ($(fa)!0.35!(fb)$);
			\coordinate (v') at ($(fa)!0.75!(fb)$);
			\coordinate (u'') at ($(fa)!0.40!(fb)$);
			\coordinate (v'') at ($(fa)!0.70!(fb)$);
			\draw[ultra thick,orange] (u'')--(v'');
			\draw[red] (fc)--(u');
			\draw[red] (fc)--(v');
			\path[name path=fafq] (fa)--(fq);
			\path[name path=u'fc] (u')--(fc);
			\path[name path=v'fc] (v')--(fc);
			\path[name path=u'v'] (u')--(v');
			\coordinate (fd) at ($(fa)!0.70!(fq)$);
			\fill (fd) circle (0.34375pt);
			\node[below right] at (fd) {$d$};
			\fill (fa) circle (0.34375pt);
			\fill (fb) circle (0.34375pt);
			\fill (fc) circle (0.34375pt);
			\fill[red] (u') circle (0.34375pt);
			\fill[red] (v') circle (0.34375pt);
			\node[below] at (fa) {$a$};
			\node[below] at  (fb) {$b$};
			\node[above] at  (fc) {$c$};
			
			\node[below] at  (u') {\textcolor{red}{$\widetilde{u}$}};
			\node[below] at  (v') {\textcolor{red}{$\widetilde{v}$}};

			\begin{scope}[scale=1,xshift=-1.2cm]

				\draw[ultra thick,orange] (u'')--(v'');
			\end{scope}
		\end{tikzpicture}\hfill
		\begin{tikzpicture}[scale=4]
			\coordinate (fa) at (0,0);
			\coordinate (fb) at (1,0);
			\coordinate (fc) at (.5,1);
			\draw (fa)--(fb)--(fc)--(fa);
			\coordinate (fq) at ($(fc)!0.35!(fb)$);
			\draw[cyan] (fa)--(fq);
			\coordinate (u') at ($(fa)!0.35!(fb)$);
			\coordinate (v') at ($(fa)!0.75!(fb)$);
			\coordinate (u'') at ($(fa)!0.40!(fb)$);
			\coordinate (v'') at ($(fa)!0.70!(fb)$);
			\draw[ultra thick,orange] (u'')--(v'');
			\draw[red] (fc)--(u');
			\draw[red] (fc)--(v');
			\draw (fc)--(u');
			\draw (fc)--(v');
			\path[name path=fafq] (fa)--(fq);
			\path[name path=u'fc] (u')--(fc);
			\path[name path=v'fc] (v')--(fc);
			\path[name path=u'v'] (u')--(v');
			\coordinate (fd) at ($(fa)!0.70!(fq)$);
			\fill (fd) circle (0.34375pt);
			\node[below right] at (fd) {$d$};
			\fill (fa) circle (0.34375pt);
			\fill (fb) circle (0.34375pt);
			\fill (fc) circle (0.34375pt);
			\fill[red] (fq) circle (0.34375pt);
			\fill[red] (u') circle (0.34375pt);
			\fill[red] (v') circle (0.34375pt);
			\fill (u') circle (0.34375pt);
			\fill (v') circle (0.34375pt);
			\node[below] at (fa) {$a$};
			\node[below] at  (fb) {$b$};
			\node[above] at  (fc) {$c$};
			\node[right,red] at  (fq) {$q$};
			
			\node[below] at  (u') {\textcolor{red}{$\widetilde{u}$}};
			\node[below] at  (v') {\textcolor{red}{$\widetilde{v}$}};
			\node[below] at  (u') {$\widetilde{u}$};
			\node[below] at  (v') {$\widetilde{v}$};

			\begin{scope}[scale=1,xshift=-1.2cm]

				\draw[blue] (fb)--(fc);
				\draw[ultra thick,orange] (u'')--(v'');
			\end{scope}
		\end{tikzpicture}\\
		\begin{tikzpicture}[scale=4]
			\coordinate (fa) at (0,0);
			\coordinate (fb) at (1,0);
			\coordinate (fc) at (.5,1);
			\draw (fa)--(fb)--(fc)--(fa);
			\coordinate (fq) at ($(fc)!0.35!(fb)$);
			\draw[cyan] (fa)--(fq);
			\draw (fa)--(fq);
			\coordinate (u') at ($(fa)!0.35!(fb)$);
			\coordinate (v') at ($(fa)!0.75!(fb)$);
			\coordinate (u'') at ($(fa)!0.40!(fb)$);
			\coordinate (v'') at ($(fa)!0.70!(fb)$);
			\draw[ultra thick,orange] (u'')--(v'');
			\draw[red] (fc)--(u');
			\draw[red] (fc)--(v');
			\draw (fc)--(u');
			\draw (fc)--(v');
			\path[name path=fafq] (fa)--(fq);
			\path[name path=u'fc] (u')--(fc);
			\path[name path=v'fc] (v')--(fc);
			\path[name path=u'v'] (u')--(v');
			\coordinate (fd) at ($(fa)!0.70!(fq)$);
			\fill (fd) circle (0.34375pt);
			\node[below right] at (fd) {$d$};
			\path[name intersections={of=fafq and u'fc,by={[label=left:\textcolor{red}{$i(\widetilde{u})$}]q1}}];
			\path[name intersections={of=fafq and v'fc,by={[label=left:\textcolor{red}{$i(\widetilde{v})$}]q2}}];
			\path[name path=u'q2] (u')--(q2);
			\path[name path=v'q1] (v')--(q1);
			\fill (fa) circle (0.34375pt);
			\fill (fb) circle (0.34375pt);
			\fill (fc) circle (0.34375pt);
			\fill[red] (fq) circle (0.34375pt);
			\fill (fq) circle (0.34375pt);
			\fill[red] (u') circle (0.34375pt);
			\fill[red] (v') circle (0.34375pt);
			\fill (u') circle (0.34375pt);
			\fill (v') circle (0.34375pt);
			\fill[red] (q1) circle (0.34375pt);
			\fill[red] (q2) circle (0.34375pt);
			\node[below] at (fa) {$a$};
			\node[below] at  (fb) {$b$};
			\node[above] at  (fc) {$c$};
			\node[right,red] at  (fq) {$q$};
			\node[right] at  (fq) {$q$};
			
			\node[below] at  (u') {\textcolor{red}{$\widetilde{u}$}};
			\node[below] at  (v') {\textcolor{red}{$\widetilde{v}$}};
			\node[below] at  (u') {$\widetilde{u}$};
			\node[below] at  (v') {$\widetilde{v}$};

			\begin{scope}[scale=1,xshift=-1.2cm]

				\draw[blue] (fb)--(fc);
				\draw[ultra thick,orange] (u'')--(v'');
			\end{scope}
		\end{tikzpicture}\hfill
        \begin{tikzpicture}[scale=4]
			\coordinate (fa) at (0,0);
			\coordinate (fb) at (1,0);
			\coordinate (fc) at (.5,1);
			\draw (fa)--(fb)--(fc)--(fa);
			\coordinate (fq) at ($(fc)!0.35!(fb)$);
			\draw[cyan] (fa)--(fq);
			\draw (fa)--(fq);
			\coordinate (u') at ($(fa)!0.35!(fb)$);
			\coordinate (v') at ($(fa)!0.75!(fb)$);
			\coordinate (u'') at ($(fa)!0.40!(fb)$);
			\coordinate (v'') at ($(fa)!0.70!(fb)$);
			\draw[ultra thick,orange] (u'')--(v'');
			\draw[red] (fc)--(u');
			\draw[red] (fc)--(v');
			\draw (fc)--(u');
			\draw (fc)--(v');
			\path[name path=fafq] (fa)--(fq);
			\path[name path=u'fc] (u')--(fc);
			\path[name path=v'fc] (v')--(fc);
			\path[name path=u'v'] (u')--(v');
			\coordinate (fd) at ($(fa)!0.70!(fq)$);
			\fill (fd) circle (0.34375pt);
			\node[below right] at (fd) {$d$};
			\path[name intersections={of=fafq and u'fc,by={[label=left:\textcolor{red}{$i(\widetilde{u})$}]q1}}];
			\path[name intersections={of=fafq and v'fc,by={[label=left:\textcolor{red}{$i(\widetilde{v})$}]q2}}];
			\path[name intersections={of=fafq and u'fc,by={[label=left:$i(\widetilde{u})$]q1}}];
			\path[name intersections={of=fafq and v'fc,by={[label=left:$i(\widetilde{v})$]q2}}];
			\draw[red] (u')--(q2);
			\draw[red] (v')--(q1);
			\path[name path=u'q2] (u')--(q2);
			\path[name path=v'q1] (v')--(q1);
			\path[name intersections={of=u'q2 and v'q1,by={[label=left:\textcolor{red}{$t$}]t}}];
			\fill (fa) circle (0.34375pt);
			\fill (fb) circle (0.34375pt);
			\fill (fc) circle (0.34375pt);
			\fill[red] (fq) circle (0.34375pt);
			\fill (fq) circle (0.34375pt);
			\fill[red] (u') circle (0.34375pt);
			\fill[red] (v') circle (0.34375pt);
			\fill (u') circle (0.34375pt);
			\fill (v') circle (0.34375pt);
			\fill[red] (t) circle (0.34375pt);
			\fill[red] (q1) circle (0.34375pt);
			\fill[red] (q2) circle (0.34375pt);
			\fill (q1) circle (0.34375pt);
			\fill (q2) circle (0.34375pt);
			\node[below] at (fa) {$a$};
			\node[below] at  (fb) {$b$};
			\node[above] at  (fc) {$c$};
			\node[right,red] at  (fq) {$q$};
			\node[right] at  (fq) {$q$};
			
			\node[below] at  (u') {\textcolor{red}{$\widetilde{u}$}};
			\node[below] at  (v') {\textcolor{red}{$\widetilde{v}$}};
			\node[below] at  (u') {$\widetilde{u}$};
			\node[below] at  (v') {$\widetilde{v}$};

			\begin{scope}[scale=1,xshift=-1.2cm]
				
				\draw[blue] (fb)--(fc);
				\draw[ultra thick,orange] (u'')--(v'');
			\end{scope}
		\end{tikzpicture}\hfill
		\begin{tikzpicture}[scale=4]
			\coordinate (fa) at (0,0);
			\coordinate (fb) at (1,0);
			\coordinate (fc) at (.5,1);
			\draw (fa)--(fb)--(fc)--(fa);
			\coordinate (fq) at ($(fc)!0.35!(fb)$);
			\draw[cyan] (fa)--(fq);
			\draw (fa)--(fq);
			\coordinate (u') at ($(fa)!0.35!(fb)$);
			\coordinate (v') at ($(fa)!0.75!(fb)$);
			\coordinate (u'') at ($(fa)!0.40!(fb)$);
			\coordinate (v'') at ($(fa)!0.70!(fb)$);
			\draw[ultra thick,orange] (u'')--(v'');
			\draw[red] (fc)--(u');
			\draw[red] (fc)--(v');
			\draw (fc)--(u');
			\draw (fc)--(v');
			\path[name path=fafq] (fa)--(fq);
			\path[name path=u'fc] (u')--(fc);
			\path[name path=v'fc] (v')--(fc);
			\path[name path=u'v'] (u')--(v');
			\coordinate (fd) at ($(fa)!0.70!(fq)$);
			\fill (fd) circle (0.34375pt);
			\node[below right] at (fd) {$d$};
			\path[name intersections={of=fafq and u'fc,by={[label=left:\textcolor{red}{$i(\widetilde{u})$}]q1}}];
			\path[name intersections={of=fafq and v'fc,by={[label=left:\textcolor{red}{$i(\widetilde{v})$}]q2}}];
			\path[name intersections={of=fafq and u'fc,by={[label=left:$i(\widetilde{u})$]q1}}];
			\path[name intersections={of=fafq and v'fc,by={[label=left:$i(\widetilde{v})$]q2}}];
			\draw[red] (u')--(q2);
			\draw[red] (v')--(q1);
			\draw (u')--(q2);
			\draw (v')--(q1);
			\path[name path=u'q2] (u')--(q2);
			\path[name path=v'q1] (v')--(q1);
			\path[name intersections={of=u'q2 and v'q1,by={[label=left:\textcolor{red}{$t$}]t}}];
			\path[name intersections={of=u'q2 and v'q1,by={[label=left:$t$]t}}];
			\fill (fa) circle (0.34375pt);
			\fill (fb) circle (0.34375pt);
			\fill (fc) circle (0.34375pt);
			\fill[red] (fq) circle (0.34375pt);
			\fill (fq) circle (0.34375pt);
			\fill[red] (u') circle (0.34375pt);
			\fill[red] (v') circle (0.34375pt);
			\fill (u') circle (0.34375pt);
			\fill (v') circle (0.34375pt);
			\fill[red] (t) circle (0.34375pt);
			\fill (t) circle (0.34375pt);
			\fill[red] (q1) circle (0.34375pt);
			\fill[red] (q2) circle (0.34375pt);
			\fill (q1) circle (0.34375pt);
			\fill (q2) circle (0.34375pt);
			\node[below] at (fa) {$a$};
			\node[below] at  (fb) {$b$};
			\node[above] at  (fc) {$c$};
			\node[right,red] at  (fq) {$q$};
			\node[right] at  (fq) {$q$};
			
			\node[below] at  (u') {\textcolor{red}{$\widetilde{u}$}};
			\node[below] at  (v') {\textcolor{red}{$\widetilde{v}$}};
			\node[below] at  (u') {$\widetilde{u}$};
			\node[below] at  (v') {$\widetilde{v}$};

			\coordinate (help) at ($(fc)!1.5!(t)$);
			\path[name path=fchelp] (fc)--(help);
			\path[name intersections={of=u'v' and fchelp,by={[label=below:\textcolor{red}{$w'$}]gg}}];
			\fill[red] (gg) circle (0.34375pt);
			\draw[cyan] (fc)--(gg);
			\begin{scope}[scale=1,xshift=-1.2cm]
				
				\draw[blue] (fb)--(fc);
				\draw[blue] (fa)--(fb);
				\draw[ultra thick,orange] (u'')--(v'');
				\fill (gg) circle (0.34375pt);
			\end{scope}
		\end{tikzpicture}
        \endgroup
        \caption{We assume by contradiction that there is a set with nonempty interior in~$[a,b]$ that does not intersect~$D$.
        It is depicted as the thick line segment in the side~$[a,b]$.
        The figures show how to obtain a contradiction by constructing a point~$w'$ in this set.}\label{FiguraJedna}
	\end{figure}
	It is easy to check that $D$ is dense in each of the sides of the triangle, therefore it is actually dense in the whole triangle. Here the fact that $d$ is in the interior of $\Delta$ is essential.
	It is clear that $f$ and $g$ are uniquely determined on $D$ (being monotone mappings), therefore also $f \rest D = g \rest D$.
	Here the density of $D$ in each side of $\Delta$ plays a crucial role, because a monotone mapping on a line segment is determined by any dense subset.
	
	Next, that $f \rest D = g \rest D$ and that $g$ restricted to the boundary of~$\Delta$ is continuous shows that $f \rest D$ has a unique monotone extension to the boundary of $\Delta$.
	Finally, $f=g$ on the boundary of $\Delta$, therefore $f = g$ on $\Delta$. A monotone mapping is totally determined on an open set, therefore $f=g$ on $G$.
\end{pf}

As a consequence, we see that the image of a monotone isomorphism of a triangle is again a triangle.
This easily extends to higher dimensions, see Theorem~\ref{THMsimplices} below.

\section{Bad monotone mappings}\label{SectZlosciBad}

One can say that monotone mappings between linearly ordered sets are relatively well understood, at least when it comes to subsets of the real line.
Concerning Euclidean betweenness, projective transformations provide the class of ``nice'' monotone maps of convex sets. In fact, these are the continuous ones, when the domain is not contained in a line. On the other hand, there are many ``bad'' monotone mappings, namely, those that are neither continuous nor one-to-one. Our aim in this section is to classify ``bad'' monotone mappings of convex planar sets.

Below are some basic examples of discontinuous monotone mappings of a triangle.

\begin{ex}\label{ExmplDwaJedne}

	Let $\Delta = [a,b,c]$ be a triangle. Define $\map f \Delta \Delta$ by setting $f(x)=x$, if $x \in [a,b]$ and $f(x) = c$, otherwise.
	It is rather clear that $f$ is monotone.

\end{ex}

\begin{ex}\label{ExmplDwwaDwwaa}

	Let $\Delta = [a,b,c]$ be a triangle.
	Fix $v \in [b, c] \setminus \{b,c\}$.
	Define $\map{g}{\Delta}{\Delta}$ by setting
	$g(x) = a$ if $x = a$; $g(x) = b$ if $x \ne a$ is in the triangle $[a,b,v]$; and $g(x) = c$, otherwise.	It is easy to see that $g$ is monotone.	

\end{ex}

In the previous example, one can take any planar convex set with at least one extreme point and with nonempty interior. Such a set admits a monotone mapping whose image consists of three non-collinear points.

\begin{ex}\label{FivePoints}
	There exists a monotone mapping $\map f {[0,1]^2}{\Err^2}$ whose image is
	\begin{equation*}
\bX = \{ (0,0), (1,0), (0,1), (1,1), (1/2, 1/2) \}.
\end{equation*}
	We define $f$ as follows:
    \begin{equation*}
	f(x_0,x_1) = \begin{cases}
	(0,1) \quad & \text{if } x_0 < x_1,\\
	(1,0) & \text{if } x_0 > x_1,\\
	(1/2,1/2) & \text{if } 0 < x_0 = x_1 < 1,\\
	(0,0) & \text{if } 0 = x_0 = x_1,\\
	(1,1) & \text{if } x_0 = x_1 = 1.
	\end{cases}
    \end{equation*}
	We check that $f$ is monotone.
	Fix $a,b \in [0,1]^2$ and $s \in [a,b]$.
	We may assume that $a \ne s \ne b$, which also rules out that $s$ is an extreme point of $[0,1]^2$.
	If $s$ is not in the diagonal $D = \setof{(x,x)}{0 \loe x \loe 1}$ then one of the points $a$, $b$ is on the same side of $D$ as $s$ (above or below) and hence $f(s) = f(a)$ or $f(s) = f(b)$.
	Finally, if $s \in D$ then, knowing that $s \ne (0,0)$ and $s \ne (1,1)$, we conclude that $f(s) = (1/2,1/2)$ and either $\img f {\dn a b} = \dn{(0,1)}{(1,0)}$ or $\img f {\dn ab} = \dn{(0,0)}{(1,1)}$.
\end{ex}

\section{A classification}

We now know that monotone mappings from planar convex sets could be pathological. On the other hand, we are able to classify them into three pairwise disjoint classes. Recall that without the requirement of convexity one can get monotone isomorphisms between very different subsets of Euclidean spaces. On the other hand, when it comes to planar sets, we may restrict attention to mappings into the real plane. Namely:

\begin{lm}\label{LMdsvsdf}
	Every monotone mapping $\map f X V$ with $X \subs \Err^2$ convex, has the property that $\img f X$ is contained in a plane.	
\end{lm}

\begin{pf}
	If all triples in $\img fX$ are collinear then $\img fX$ is contained in a line. Otherwise, fix $a,b,c \in X$ such that $f(a), f(b), f(c)$ are not collinear.
	Fix $x \in X$. Let $L$ be a line passing through $x$ and intersecting the triangle $[a,b,c]$ in more than one point and passing through one of the vertices $a,b,c$. So there are two distinct points $u,v$ such that $u \in \{a,b,c\}$, $v$ is on the side opposite to $u$, and $u,v,x$ are collinear.
	Then $f(u) \ne f(v)$ and $f(x)$ is on the line passing through $f(u), f(v)$. This line is contained in the affine plane $P$ spanned by $f(a),f(b),f(c)$. This shows that $\img fX \subs P$.
\end{pf}

Below is the announced classification.
Recall that a set $S$ is \emph{convex independent} if $x \notin \conv(S \setminus \sn x)$ for every $x \in S$.

\begin{tw}\label{THMgdyifgi}
	Assume $G \subs \Err^2$ is a nonempty convex set and $\map f G {\Err^2}$ is a monotone mapping.
	Then exactly one of the following possibilities occur.
	\begin{enumerate}[itemsep=0pt]
		\item[{\rm(1)}] $\img f G \subs L \cup \sn q$, where $L \subs \Err^2$ is a line and $q \in \Err^2$.
		\item[{\rm(2)}] $\Int fG \nnempty$ and $f$ is a partial homography.
		\item[{\rm(3)}] $\img fG = \{a,b,c,d,e\}$, where $a,b,c,d$ are convex independent and $\sn e = [a,c] \cap [b,d]$.
	\end{enumerate}
\end{tw}

\begin{pf}
	Obviously, (1) -- (3) are mutually exclusive. The second part of (2) follows from Theorem~\ref{THMhomographies}.
	
	Suppose $f$ satisfies neither (1) nor (2).
	By Theorem~\ref{THMhomographies},
	for every $x,y,z \in G$ either the points $f(x), f(y), f(z)$ are collinear or else
    \begin{equation*}
      \img f G \cap \Int [f(x), f(y), f(z)] = \emptyset.
    \end{equation*}
	We claim that $\img f G$ contains four convex independent points $a,b,c,d$.
	Suppose otherwise.
	Since $\img f G$ is not contained in a single line, there are non-collinear points $x,y,z \in \img f G$.
	Let $A = L_{x,y} \cup L_{x,z} \cup L_{y,z}$, where $L_{s,t}$ denotes the line passing through the points $s$~and~$t$.
	If there is $u \in \img f G \setminus A$ then the set $\{x,y,z,u\}$ is convex independent, because by the remarks above, none of these points is in the interior of the triangle spanned by the other three, and we have found our pair of points.
	Now we assume $\img f G \subs A$.
	Again using the fact that $f$ does not satisfy (1), we find $u,v \in \img f G \setminus \{x,y,z\}$ such that $u$ and $v$ are on different lines forming the set $A$.
	Permuting $x,y,z$ if necessary, we may assume that $u \in L_{x,y}$ and $v \in L_{x,z}$.
	Checking all possible configurations, we conclude that there are always four convex independent points.
	Indeed, if $x \in [u,y]$ then necessarily $x \in [v,z]$ since otherwise either $v \in \Int [u,y,z]$ or $z \in \Int [u,y,v]$.
	Thus, in this case the set $\{u,v,y,z\}$ is convex independent.
	We get the same conclusion if $x \in [v,z]$ (this is symmetric).
	In the remaining cases we conclude as well that the set $\{u,v,y,z\}$ is convex independent.
	
	We have proved that there is a convex independent set $\{a,b,c,d\} \subs \img f G$ consisting of four elements.
	Let $\{a',b',c',d'\} \subs G$ be such that $a = f(a')$, $b = f(b')$, $c = f(c')$, and $d = f(d')$.
	Then the four-element set $\{a',b',c',d'\}$ is convex independent, therefore, possibly permuting the points, we may assume that $[a',c'] \cap [b',d'] = \sn {e'}$.
	Let $e = f(e')$.
	We claim that $\img f G = \{a,b,c,d,e\}$.
	
	Clearly, $\sn e = [a,c] \cap [b,d]$.
	Assume by contradiction that $w \in \img f G \setminus \{a,b,c,d,e\}$.
	If $w$ is in the interior of the angle $\measuredangle(e,a,b)$ then $e$ is in the interior of the triangle $[c,d,w]$, a contradiction.
	Similarly, we show that necessarily $w \in L_{a,c} \cup L_{b,d}$.
	However, in this case we find, checking all possible configurations (by symmetry, it suffices to check only the cases where $w\in [a,e]$ or $a\in [w,e]$) that the set $\{a,b,c,d,w\}$ contains a point that is in the interior of the triangle spanned by other three points, again obtaining a contradiction.
	This completes the proof.
\end{pf}

The result above holds also when the range of the mapping $f$ is an arbitrary real vector space, because due to Lemma~\ref{LMdsvsdf} the image of $G$ is contained in an affine subplane. On the other hand, property~(2) would need to be rephrased, namely, the interior would need to be computed in the plane spanned by $\img fG$. Actually, property~(2) is equivalent to the existence of four non-collinear points such that one of them is in the interior of the triangle spanned by the other three points.

In case that $G$ is open, we have just a dichotomy.

\begin{wn}
	Assume $G \subs \Err^2$ is a nonempty open convex set and $\map f G {\Err^2}$ is monotone. Then either $\img fG$ is contained in a line or else $f$ is a partial homography.
\end{wn}

\begin{pf}
	Suppose (1) or (3) holds and $\img fG$ is not contained in a line. Then $\img fG$ can be decomposed into three nonempty half-spaces (at least one of them being a singleton). Their pre-images form a partition of $G$ into three nonempty half-spaces, which is impossible, because $G$ is open.
\end{pf}

\begin{wn}\label{WNdnfosdogi}
	Every one-to-one monotone mapping of a convex planar set into~$\Err^2$ is a partial homography.
\end{wn}

\begin{pf}
	Let $\map{f}{G}{\Err^2}$ be monotone and one-to-one. By Theorem~\ref{THMsimplices}, $f$ restricted to any triangle contained in $G$ is a partial homography. Hence $f$ is a partial homo\-gra\-phy, because it is determined by its restriction to an arbitrarily small convex open subset.
\end{pf}

\ignoreSpellCheck{Hou} and \ignoreSpellCheck{McColm}~\cite{HM08} showed that every monotone affine self-mapping of the plane is affine. This follows from the above corollary, because a partial homography whose domain is the whole plane is necessarily affine.

\section{One-to-one monotone mappings}

In this section we address the question which planar sets admit a monotone one-to-one mapping into the real line.
The first result exhibits a natural influence of the cardinality continuum.

\begin{tw}
	Let $S$ be a subset of a real vector space $V$ and let $L$ be a $1$\nobreakdash-dimensional linear subspace of $V$. If the cardinality of $S$ is strictly less than the continuum then there exists a linear map $\map P V L$ such that $P \rest S$ is one-to-one.
	In particular, $S$ admits a one-to-one monotone mapping into the real line.
\end{tw}

\begin{pf}
	By Proposition~\ref{propreducn2}, we may assume that $V$ is $2$\nobreakdash-dimensional.
	Let $\El$ be the family of lines passing through at least two points of $S$. Then $\El$ has cardinality strictly less than the continuum, therefore there is a line $K$ passing through the origin, not parallel to any line from the family $\El$ and different from $L$. Let $P$ be the projection onto $L$ along $K$. Then $P\rest S$ is obviously one-to-one.
\end{pf}

The next result says that convex sets never have the above property, unless they are contained in a single line.

\begin{tw}\label{THMtrojkatNieProsty}
	Let $G \subs \Err^2$ be a  set with nonempty interior, let $\map f G Y$ be a one-to-one monotone mapping, where $Y$ is a linearly ordered set.
	Then $Y$ contains continuum many pairwise disjoint intervals, each of them of cardinality at least the continuum.
\end{tw}

\begin{pf}
	We may assume $G$ is open and convex, replacing it by a small enough open ball. We may also assume that $Y = \img f G$.
	Given $t \in Y$, define
    \begin{equation*}
	  H_t = \inv{f}{(-\infty, t]} = \inv{f}{\setof{y \in Y}{y \loe t}}.
    \end{equation*}
	Then $H_t$ is a nontrivial (proper, infinite) half-space in $G$, therefore it is determined by a unique line $L_t$. Specifically, $H_t = U_t \cup I_t$, where $U_t$ is an open half-space disjoint from $L_t$ and $I_t \subs L_t \cap G$ is an interval.
		Note that $I_t \nnempty$, because, denoting by $s$ the unique point such that $f(s) = t$, we have that $H_t \setminus \sn s = \inv f {(-\infty, t)}$ is a half-space, too, therefore it contains $U_t$.
	Now let
    \begin{equation*}
      \El = \setof{L_t}{t \in Y}.
    \end{equation*}
	Note that $G \subs \bigcup \El$.
	Indeed, fix $x \in G$ and let $t = f(x)$.
	Then $L_t \in \El$ and $x \in L_t$.
	It follows that $\El$ has cardinality continuum.
	
	Now suppose that $L_t \cap G = (c_t,d_t)$, where $f(c_t) < f(d_t)$. Let $J_t = \img f {(c_t,d_t)}$.
	We claim that if $s < t$ are such that $L_s \ne L_t$ then $J_s \cap J_t = \emptyset$.
	Indeed, in this case $(c_t, d_t)$ it is disjoint from $H_s$, therefore $f(d_s) < f(c_t)$.
	Finally, $\setof{J_t}{t \in Y}$ is a family of pairwise disjoint open intervals in $Y$, and it has cardinality $|\El|$, which is the continuum.	
\end{pf}

\begin{tw}\label{THMsimplices}
	Let $\Delta_m$ be an $m$\nobreakdash-dimensional simplex, where $m \goe 2$.
	Let $\map f {\Delta_m} V$ be a one-to-one monotone mapping, where $V$ is a real vector space.
	Then $\img f {\Delta_m}$ is an $m$\nobreakdash-dimensional simplex and $f$ is a partial homography.
\end{tw}

\begin{pf}
	First of all, note that $f$ is an isomorphism onto its image.
	Indeed, if $x,y,z$ are not collinear and $f(x), f(y), f(z)$ are collinear, then $f$ restricted to the triangle $[x,y,z]$ would be a one-to-one monotone mapping into a line segment, which contradicts Theorem~\ref{THMtrojkatNieProsty}.
	
	Let us start with the case $m=2$. Replacing $f$ by the composition with a suitable linear transformation, we may assume that $\img f{\Delta_2} \subs \Err^2$.	
	Now, since we know that $f$ is an isomorphism, the assumptions of Theorem~\ref{THMhomographies} are fulfilled, therefore $f$ is a partial homography.
	In particular, its image is the triangle spanned by the images of the vertices of $\Delta_2$.
	Clearly, the same conclusion is true if $\Err^2$ is replaced by any other (not necessarily $2$\nobreakdash-dimensional) real vector space.
	
	Now assume $m>2$.
    Recall that the \emph{Helly number} is the minimal $k$ such that every finite family of convex sets with empty intersection contains at most $k$ sets with empty intersection.
    A classical result of Helly says that the Helly number of $\Err^m$ is $m+1$.
By the inductive hypothesis, $\img f F$ is an $(m-1)$\nobreakdash-dimensional simplex for each $(m-1)$\nobreakdash-dimensional face $F \subs \Delta_m$.
	Thus, supposing $\img{f}{\Delta_m}$ is contained in an $(m-1)$\nobreakdash-dimensional space, the Helly number contradicts the fact that $f$ is a monotone isomorphism. Namely, all $(m-1)$\nobreakdash-dimensional faces of $\Delta_m$ have empty intersection, while their images have nonempty intersection, because every $m$ of them have nonempty intersection.
\end{pf}

Note that a planar convex set with empty interior consists of collinear points, therefore the theorem above says precisely which planar convex sets admit a one-to-one monotone mapping into $\Err$.

Concerning arbitrary planar sets, we do not know any characterization (it might be the case that no reasonable one exists). Let us first look at sets consisting of (or containing) finitely (or even countably) many lines.

\begin{prop}
	Assume $S \subs \Err^2$ consists of countably many pairwise parallel lines. Then $S$ admits a one-to-one monotone mapping into the real line.
\end{prop}

\begin{pf}
	Assume $S = \bigcup_{i \in J}L_i$, where each $L_i$ is a horizontal line and $J$ is countable. We may assume that $J \subs \Err$ and $L_i = (\Err \times \sn 0) + (0,i)$ for $i \in J$.
    Given two subsets $A$~and~$B$ of the real line, we write $A<B$ if for all $a\in A$ and $b\in B$, we have $a<b$.
	Let $\sett{U_i}{i \in J}$ be a family of pairwise disjoint open intervals in $\Err$ satisfying $U_i < U_j$ whenever $i<j$, $i,j \in J$.
	This is possible, because $J$ is countable. More precisely, each $U_i$ may be chosen among open intervals removed during the construction of the standard Cantor set, since these intervals are ordered as the rational numbers.
	
	Now let $\map f S \Err$ be such that $f \rest L_i$ is an increasing mapping onto $U_i$ in the sense that if $x<y$ and if $(x,i)$~and~$(y,i)$ are in~$L_i$, then $f(x,i)\leq f(y,i)$. Then $f$ is monotone, because if $a,b,c$ are such that $B(a,c,b)$ and $a \in L_i$, $c \in L_j$, $b \in L_k$ then either $i \loe j \loe k$ or $k \loe j \loe i$.
\end{pf}

An interval is \emph{non-degenerate} if it contains at least two points.
\begin{tw}\label{threeLinesToLine}
	Let $L_0, L_1, L_2$ be different lines in the plane $\Err^2$ such that each two of them intersect,
	let $S = (L_0 \cup L_1 \cup L_2) \setminus B$, where $B$ is a bounded set. Assume $Y$ is a linearly ordered set and $\map f S Y$ is monotone and one-to-one.
	Then $Y$ contains continuum many pairwise disjoint non-degenerate closed intervals.
	In particular, $S$ does not admit any one-to-one monotone mapping into $\Err$.
\end{tw}

\begin{pf}
\begin{figure}[!b]
\centering
    \begingroup
    \iffarben
    \else
    \colorlet{orange}{black}
    \colorlet{red}{black}
    \colorlet{cyan}{black}
    \colorlet{blue}{black}
    \colorlet{green}{black}
    \fi
\begin{subfigure}{0.48\textwidth}
\begin{tikzpicture}[scale=.40]

\node (v1) at (-7.5,-0.5) {};
\node (v2) at (5.5,-0.5) {};
\node (v4) at (4,5) {};
\node (v3) at (-9.5,-7.75) {};
\node (v5) at (-5.5,2.5) {};
\node (v6) at (6.5,-5) {};
\draw  (v3) edge (v4);
\draw  (v5) edge (v6);

\node (v7) at (-2.5,0) {};
\node (v10) at (-2.5,-1.5) {};
\node (v9) at (-0.5,-1.5) {};
\node (v8) at (-0.5,0) {};

\draw  (v7.center) edge (v8.center);
\draw  (v8.center) edge (v9.center);
\draw  (v9.center) edge (v10.center);
\draw  (v10.center) edge (v7.center);
\node (v11) at (-5,-0.5) {};
\node[below] (v12) at (-5,-0.5) {$x$};;
\draw[red,ultra thick] (v11.center) edge (v1.center);
\node (li) at (-2.5,-0.5) {};
\draw[green,thick,densely dashed] (v11.center) edge (li.center);
\node (ri) at (-0.5,-0.5) {};
\draw[green,thick,densely dashed] (v2.center) edge (ri.center);
\node[below right] at (-4,-2.5) {$L_2^-$};
\node[below right] at (-7.5,-0.5) {\!\!\!\!$L_1^-$};
\node[below] at (-4.5,2)  {$L_0^-$};
\node[below] at (2,3) {$L_2^+$};
\node[below] at (4,-0.5) {$L_1^+$};
\node[below] at (5,-4) {$L_0^+$};
\draw[orange,dotted,thick] (3,-5.5) edge (-9,2);
\draw[orange,dotted,thick] (-9.5,-4.75) edge (4,8);
\fill[gray]  (v7) rectangle (v9);
\node [blue,regular polygon, regular polygon sides=5,scale=0.3,fill] at (v11) {};
\end{tikzpicture}
\caption{The line \tikz[baseline=-.6ex]{\draw[red,ultra thick] (0,0)--(1.0,0)} is a subset of the points mapped to values smaller or equal than~$f(x)$, while the set \tikz[baseline=-.6ex]{\draw[green,thick,densely dashed] (0,0)--(1.0,0)} is mapped to values larger than~$f(x)$. Finally, the lines \tikz[baseline=-.6ex]{\draw[orange,dotted,thick] (0,0)--(1.0,0)} are the lines through~$x$ that are parallel to $L_0$~and~$L_2$.}
\end{subfigure}\hfill
\begin{subfigure}{0.48\textwidth}
\begin{tikzpicture}[scale=.40]

\node (v1) at (-7.5,-0.5) {};
\node (v2) at (5.5,-0.5) {};
\node (v4) at (4,5) {};
\node (v3) at (-9.5,-7.75) {};
\node (v5) at (-5.5,2.5) {};
\node (v6) at (6.5,-5) {};
\draw (v3) -- (v4);
\draw  (v5) edge (v6);

\node (v7) at (-2.5,0) {};
\node (v10) at (-2.5,-1.5) {};
\node (v9) at (-0.5,-1.5) {};
\node (v8) at (-0.5,0) {};
\draw  (v7.center) edge (v8.center);
\draw  (v8.center) edge (v9.center);
\draw  (v9.center) edge (v10.center);
\draw  (v10.center) edge (v7.center);
\node (v11) at (-5,-0.5) {};
\node[below] (v12) at (-5,-0.5) {$x$};
\draw[red,ultra thick] (v11.center) edge (v1.center);
\node (li) at (-2.5,-0.5) {};
\draw[green,thick,densely dashed] (v11.center) edge (li.center);
\node (ri) at (-0.5,-0.5) {};
\draw[green,thick,densely dashed] (v2.center) edge (ri.center);
\node[below right] at (-4.3,-2.5) {\!\!\!\!$L_2^-$};
\node[below right] at (-7.5,-0.5) {\!\!\!\!$L_1^-$};
\node[below] at (-4.5,2)  {$L_0^-$};
\node[below] at (2,3) {$L_2^+$};
\node[below] at (4,-0.5) {$L_1^+$};
\node[below] at (5,-4) {$L_0^+$};
\draw[orange,dotted,thick] (3,-5.5) -- (-9,2);
\draw[orange,dotted,thick] (-9.5,-4.75) edge (4,8);
\path[name path=p2] (3,-5.5) --(-9,2);
\path[name path=34] (v3) --(v4);
\path[name intersections={of=p2 and 34,by={A}}];
\node (B) at (0,2) {};
\node (C) at (0,-4) {};
\draw[red] let \p{A}=(A), \p{B}=(B) in (A)--(\x{A},\y{B});
\draw[red] let \p{A}=(A), \p{C}=(C) in (A)--(\x{A},\y{C});
\fill[gray]  (v7) rectangle (v9);
\node [blue,regular polygon, regular polygon sides=5,scale=0.3,fill] at (v11) {};
\end{tikzpicture}
\caption{If we assume that the intersection points of the vertical line with $L_0$~and~$L_2$ are mapped to values smaller than~$f(x)$, then also the intersection point with $L_1$, a contradiction.\\ $ $}
\end{subfigure}
\begin{subfigure}{0.48\textwidth}
\begin{tikzpicture}[scale=.40]

\node (v1) at (-7.5,-0.5) {};
\node (v2) at (5.5,-0.5) {};
\node (v4) at (4,5) {};
\node (v3) at (-9.5,-7.75) {};
\node (v5) at (-5.5,2.5) {};
\node (v6) at (6.5,-5) {};
\draw  (v5) edge (v6);

\node (v7) at (-2.5,0) {};
\node (v10) at (-2.5,-1.5) {};
\node (v9) at (-0.5,-1.5) {};
\node (v8) at (-0.5,0) {};
\draw  (v7.center) edge (v8.center);
\draw  (v8.center) edge (v9.center);
\draw  (v9.center) edge (v10.center);
\draw  (v10.center) edge (v7.center);
\node (v11) at (-5,-0.5) {};
\node[below] (v12) at (-5,-0.5) {$x$};
\draw[red,ultra thick] (v11.center) edge (v1.center);
\node (li) at (-2.5,-0.5) {};
\draw[green,thick,densely dashed] (v11.center) edge (li.center);
\node (ri) at (-0.5,-0.5) {};
\draw[green,thick,densely dashed] (v2.center) edge (ri.center);
\node[below right] at (-5.5,-3.5) {$L_2^-$};
\node[below right] at (-7.5,-0.5) {\!\!\!\!$L_1^-$};
\node[below] at (-4.5,2)  {$L_0^-$};
\node[below] at (2,3) {$L_2^+$};
\node[below] at (4,-0.5) {$L_1^+$};
\node[below] at (5,-4) {$L_0^+$};
\draw[orange,dotted,thick] (3,-5.5) edge (-9,2);
\draw[orange,dotted,thick] (-9.5,-4.75) edge (4,8);
\path[name path=p2] (3,-5.5) --(-9,2);
\path[name path=34] (v3) --(v4);
\path[name intersections={of=p2 and 34,by={A}}];
\draw[green,thick,densely dashed] (A)--(v4);
\draw (v3)--(A);
\node[below,rotate=325] at (A) {\,\,\,\,\,\,\,\,\,\,\,\,$p_2(x)$};
\fill[gray]  (v7) rectangle (v9);
\node [blue,regular polygon, regular polygon sides=5,scale=0.3,fill] at (v11) {};
\node [blue,regular polygon, regular polygon sides=5,scale=0.3,fill] at (A) {};
\end{tikzpicture}
\caption{Without loss of generality, the part on~$L_2^-$ on the right of~$p_2(x)$ is mapped to values larger than~$f(x)$.}
\end{subfigure}\hfill
\begin{subfigure}{0.48\textwidth}
\begin{tikzpicture}[scale=.40]

\node (v1) at (-7.5,-0.5) {};
\node (v2) at (5.5,-0.5) {};
\node (v4) at (4,5) {};
\node (v3) at (-9.5,-7.75) {};
\node (v5) at (-5.5,2.5) {};
\node (v6) at (6.5,-5) {};
\draw  (v5) edge (v6);

\node (v7) at (-2.5,0) {};
\node (v10) at (-2.5,-1.5) {};
\node (v9) at (-0.5,-1.5) {};
\node (v8) at (-0.5,0) {};
\draw  (v7.center) edge (v8.center);
\draw  (v8.center) edge (v9.center);
\draw  (v9.center) edge (v10.center);
\draw  (v10.center) edge (v7.center);
\node (v11) at (-5,-0.5) {};
\node[below] (v12) at (-5,-0.5) {$x$};
\draw[red,ultra thick] (v11.center) edge (v1.center);
\node (li) at (-2.5,-0.5) {};
\draw[green,thick,densely dashed] (v11.center) edge (li.center);
\node (ri) at (-0.5,-0.5) {};
\draw[green,thick,densely dashed] (v2.center) edge (ri.center);
\node[below right] at (-4,-2.5) {$L_2^-$};
\node[below right] at (-7.5,-0.5) {\!\!\!\!$L_1^-$};
\node[below] at (-4.5,2)  {$L_0^-$};
\node[below] at (2,3) {$L_2^+$};
\node[below] at (4,-0.5) {$L_1^+$};
\node[below] at (5,-4) {$L_0^+$};
\draw[orange,dotted,thick] (3,-5.5) edge (-9,2);
\draw[orange,dotted,thick] (-9.5,-4.75) edge (4,8);
\path[name path=p2] (3,-5.5) --(-9,2);
\path[name path=34] (v3) --(v4);
\path[name intersections={of=p2 and 34,by={A}}];
\draw[green,thick,densely dashed] (A)--(v4);
\draw (v3)--(A);
\node (v13) at (-5.5,-1.7) {};
\node (v14) at (1,1.7) {};
\draw  (v13) edge (v14);
\fill[gray]  (v7) rectangle (v9);
\node [blue,regular polygon, regular polygon sides=5,scale=0.3,fill] at (v11) {};
\end{tikzpicture}
\caption{We obtain more information about the values of the images of some points on~$L_0$.\\ \hphantom{We obtain more information about the}}
\end{subfigure}
\caption{Following the proof of~Theorem~\ref{threeLinesToLine} (continued on the next page)}\label{The line}
\endgroup
\end{figure}
\begin{figure}[ht]\ContinuedFloat
    \begingroup
    \iffarben
    \else
    \colorlet{orange}{black}
    \colorlet{red}{black}
    \colorlet{cyan}{black}
    \colorlet{blue}{black}
    \colorlet{green}{black}
    \fi
\begin{subfigure}{0.48\textwidth}
\begin{tikzpicture}[scale=.40]

\node (v1) at (-7.5,-0.5) {};
\node (v2) at (5.5,-0.5) {};
\node (v4) at (4,5) {};
\node (v3) at (-9.5,-7.75) {};
\node (v5) at (-5.5,2.5) {};
\node (v6) at (6.5,-5) {};

\node (v7) at (-2.5,0) {};
\node (v10) at (-2.5,-1.5) {};
\node (v9) at (-0.5,-1.5) {};
\node (v8) at (-0.5,0) {};
\draw  (v7.center) edge (v8.center);
\draw  (v8.center) edge (v9.center);
\draw  (v9.center) edge (v10.center);
\draw  (v10.center) edge (v7.center);
\node (v11) at (-5,-0.5) {};
\node[below] (v12) at (-5,-0.5) {$x$};
\draw[red,ultra thick] (v11.center) edge (v1.center);
\node (li) at (-2.5,-0.5) {};
\draw[green,thick,densely dashed] (v11.center) edge (li.center);
\node (ri) at (-0.5,-0.5) {};
\draw[green,thick,densely dashed] (v2.center) edge (ri.center);
\node[below right] at (-4,-2.5) {$L_2^-$};
\node[below right] at (-7.5,-0.5) {\!\!\!\!$L_1^-$};
\node[below] at (-4.5,2)  {$L_0^-$};
\node[below] at (2,3) {$L_2^+$};
\node[below] at (4,-0.5) {$L_1^+$};
\node[below] at (5,-4) {$L_0^+$};
\draw[orange,dotted,thick] (3,-5.5) edge (-9,2);
\draw[orange,dotted,thick] (-9.5,-4.75) edge (4,8);
\path[name path=p2] (3,-5.5) --(-9,2);
\path[name path=34] (v3) --(v4);
\path[name intersections={of=p2 and 34,by={A}}];
\draw[green,thick,densely dashed] (A)--(v4);

\path[name path=56] (v5) --(v6);
\path[name path=or] (-9.5,-4.75)--(4,8);
\path[name intersections={of=or and 56,by={B}}];
\draw[green,thick,densely dashed] (v6)--(B);
\draw (v5)--(B);
\draw (v3)--(A);

\fill[gray]  (v7) rectangle (v9);
\node [blue,regular polygon, regular polygon sides=5,scale=0.3,fill] at (v11) {};

\end{tikzpicture}
\caption{We update the figure according to our amassed knowledge.\newline\hphantom{We update the figure according to our} \newline\hphantom{We update the figure according to our}\newline\hphantom{We update the figure according to our}}
\end{subfigure}\hfill
\begin{subfigure}{0.48\textwidth}
\begin{tikzpicture}[scale=.40]
\node (v1) at (-7.5,-0.5) {};
\node (v2) at (5.5,-0.5) {};
\node (v4) at (4,5) {};
\node (v3) at (-9.5,-7.75) {};
\node (v5) at (-5.5,2.5) {};
\node (v6) at (6.5,-5) {};

\node (v7) at (-2.5,0) {};
\node (v10) at (-2.5,-1.5) {};
\node (v9) at (-0.5,-1.5) {};
\node (v8) at (-0.5,0) {};
\draw  (v7.center) edge (v8.center);
\draw  (v8.center) edge (v9.center);
\draw  (v9.center) edge (v10.center);
\draw  (v10.center) edge (v7.center);
\node (v11) at (-5,-0.5) {};
\node (v12) at (-5,-0.5) {};
\node[below] at (-5,-0.5) {\,\,$x$};
\node (x) at (v12) {};
\draw[red,ultra thick] (v11.center) edge (v1.center);
\node (li) at (-2.5,-0.5) {};
\draw[green,thick,densely dashed] (v11.center) edge (li.center);
\node (ri) at (-0.5,-0.5) {};
\draw[green,thick,densely dashed] (v2.center) edge (ri.center);
\node[below right] at (-4,-2.5) {$L_2^-$};
\node[below right] at (-7.5,-0.5) {\!\!\!\!$L_1^-$};
\node[below] at (-2.0,2)  {\,$L_0^-$};
\node[below] at (2,3) {$L_2^+$};
\node[below] at (4,-0.5) {$L_1^+$};
\node[below] at (5,-4) {$L_0^+$};
\draw[orange,dotted,thick] (3,-5.5) edge (-9,2);
\draw[orange,dotted,thick] (-9.5,-4.75) edge (4,8);
\path[name path=p2] (3,-5.5) --(-9,2);
\path[name path=34] (v3.center) --(v4.center);
\path[name intersections={of=p2 and 34,by={A}}];

\path[name path=56] (v5) --(v6);
\path[name path=or] (-9.5,-4.75)--(4,8);
\path[name intersections={of=or and 56,by={B}}];

\node (Q) at ($(v5)!.06!(v6)$) {};
\node[above] at ($(v5)!.06!(v6)$) {$Q$};
\node (Qp) at ($(v5)!.13!(v6)$) {};
\node[above] at ($(v5)!.13!(v6)$) {$Q'$};
\path[name path=uuu] ($(Q)+(Q)-(x)$)--($(x)+(x)-(Q)+(x)-(Q)$);
\path[name path=uuuu] ($(Qp.center)+(Qp.center)-(x.center)$)--($(x.center)+(x.center)-(Qp.center)+(x.center)-(Qp.center)+(x.center)-(Qp.center)+(x.center)-(Qp.center)$);
\draw[purple,name intersections={of=34 and uuu,by={D}}];
\node[below]  at (D) {$D$};
\draw[purple,name intersections={of=34 and uuuu,by={Dp}}];
\node[below]  at (Dp) {$D'$};

\draw[green,thick,densely dashed] (Q.center) edge (v6);%
\draw[red,ultra thick] (Q.center) edge (v5);
\draw (Q.center) edge (D.center);
\draw (Qp.center) edge (Dp.center);
\draw[green,thick,densely dashed] (v4.center) edge (Dp.center);%
\draw (v3)--(Dp);
\fill[gray]  (v7) rectangle (v9);
\node [blue,regular polygon, regular polygon sides=5,scale=0.3,fill] at (v11) {};

\node [blue,regular polygon, regular polygon sides=5,scale=0.3,fill] at (D) {};
\node [blue,regular polygon, regular polygon sides=5,scale=0.3,fill] at (Dp) {};
\node [blue,regular polygon, regular polygon sides=5,scale=0.3,fill] at (Q) {};
\node [blue,regular polygon, regular polygon sides=5,scale=0.3,fill] at (Qp) {};
\end{tikzpicture}
\caption{Let $Q$ on~$L_0^-$ such that $f(Q)<f(x)$, but not much smaller.
We see that the points between $D$~and~$D'$ cannot be mapped to values larger than~$f(x)$ as this would imply that also $f(x)$ has to be larger than~$f(x)$.}
\end{subfigure}
\caption{(continued) Following the proof of~Theorem~\ref{threeLinesToLine} }
\endgroup
\end{figure}
	We may assume that $B$ is convex, closed, and contains all the intersections of our three lines (replacing $B$ by a sufficiently big closed rectangle or a closed ball).
	We may also assume that $Y = \img f S$.
	Note that each $L_i$ is partitioned into two half-lines $L_i^-$, $L_i^+$, where $\img f{L_i^-} < \img f {L_i^+}$. This is because $f$ is monotone therefore either preserves or reverses the fixed natural ordering of $L_i \setminus B$.
	
	Rotating the plane and possibly re-enumerating the lines, we may assume that $L_1$ is a horizontal line, $L_0^-$ is above $L_1$ and $L_2^-$ is below $L_1$. By this way, $L_0^+$ is below $L_1$ and $L_2^+$ is above $L_1$.
	It is easy to check that other configurations (assuming $L_1$ is horizontal) would violate the fact $f$ is a monotone injection.
	We may also assume that neither $L_0$ nor $L_2$ is vertical.
	By this way, each of our three lines has a natural \emph{horizontal} linear ordering (namely, the one induced by the first coordinate) and $f$ preserves this ordering.
	
	Given $x \in L_1^-$, denote by $p_2(x)$ the unique point in the intersection of $L_2$ with the line passing through $x$ parallel to $L_0$.
	Note that if $x$ is far enough from $B$ then $p_2(x) \in L_2^-$.
	Let us define $p_0(x) \in L_0$ in the same manner.
	Let us call $x$ \emph{relevant} if both $p_2(x) \in L_2^-$ and $p_0(x) \in L_0^-$.
	Clearly, all points of $L_1^-$ that are far enough from~$B$ are relevant.
	
	Fix a relevant point $x \in L_1^-$.
	We claim that $f(y) < f(x)$ for every $y \in L_2^-$ above $p_2(x)$ in the horizontal ordering.
	Suppose otherwise, and fix a witness $y$.
    The reader may want to have a look at Figure~\ref{The line}.
	Let $z \in L_0^+$ be such that $x, y, z$ are collinear.
	Then $f(z) > f(y) > f(x)$.
	Now take $y'$ strictly between $p_2(x)$ and $y$ and let $z' \in L_0^+$ be such that $x, y', z'$ are collinear.
	Then $f(z') < f(z)$, although it should be the opposite.
	We have used the fact $y$ was close enough to $p_2(x)$, so that $z \notin B$, but this is not a problem, as one can always replace $y$ by a point closer to $p_2(x)$. We have also used the fact that $S \cap L_2$ is open.
	
	Now fix a relevant point $x \in L_1^-$.
	Let
    \begin{equation*}
	  H_x = \setof{s \in S}{f(s) \loe f(x)}.
    \end{equation*}
	Then $H_x$ is a half-space in $S$ and by the remarks above, it is determined by a line $K_x$ either parallel to $L_0$ or $L_2$ (and then passing through either $p_0(x)$ or $p_2(x)$) or intersecting both $L_0^-$ and $L_2^-$.
	In particular, $K_x$ is not horizontal and $H_x$ contains all points of $S$ that are ``on the left side'' of the line $K_x$.
	
	Let us call a relevant point $x \in L_1^-$ \emph{tame} if $K_x$ is parallel to $L_0$ or $L_2$, otherwise let us call it \emph{wild}.
	
	Suppose that among the relevant points of $L_1^-$ continuum many of them are tame.
	Then there is a set $T \subs L_1^-$ of cardinality continuum such that for every $x \in T$ the line $K_x$ passes through $p_i(x)$ where $i\in \dn 02$ is fixed.
	Then $\sett{[f(x), f(p_i(x))]}{x \in T}$ is a family of pairwise disjoint non-degenerate closed intervals of $Y$.
	
	Now suppose that there is a set $W \subs L_1^-$ of cardinality continuum, consisting of wild points.
	If there are distinct $x, x' \in W$ such that $K_x, K_{x'}$ passes through the same point of some $L_i$ ($i \in \dn02$) then we proceed as follows.
	Reverting the picture, if necessary, assume $i=0$.
	Let $J$ be the interval between $x$ and $x'$.
	Then $\sett{[y, p_2(y)]}{y\in J}$ consists of pairwise disjoint non-degenerate closed intervals.
	
	Finally, let us assume that the lines $K_x$, $x \in W$, do not intersect in $S$.
    We denote by $q_0(x)$ its intersection with~$L_0$ and by $q_2(x)$ its intersection with $L_2$.
	Then $\sett{[f(q_0(x)), f(q_2(x))]}{x \in W}$ is a pairwise disjoint family of closed non-degenerate intervals, see Figure~\ref{position of lines}. In the last case, we have that $f(x) \in [f(q_0(x)), f(q_2(x))]$, therefore these intervals have at least three distinct points.
\begin{figure}
\centering
    \begingroup
    \iffarben
    \else
    \colorlet{orange}{black}
    \colorlet{red}{black}
    \colorlet{cyan}{black}
    \colorlet{blue}{black}
    \colorlet{green}{black}
    \fi
\begin{tikzpicture}[scale=.70]
\node (v1) at (-7.5,-0.5) {};
\node (v2) at (5.5,-0.5) {};
\node (v4) at (4,5) {};
\node (v3) at (-9.5,-7.75) {};
\node (v5) at (-5.5,2.5) {};
\node (v6) at (6.5,-5) {};

\node (v7) at (-2.5,0) {};
\node (v10) at (-2.5,-1.5) {};
\node (v9) at (-0.5,-1.5) {};
\node (v8) at (-0.5,0) {};
\draw  (v7.center) edge (v8.center);
\draw  (v8.center) edge (v9.center);
\draw  (v9.center) edge (v10.center);
\draw  (v10.center) edge (v7.center);
\node (v11) at (-5,-0.5) {};
\node (v12) at (-5,-0.5) {};
\node[below] at (-5,-0.5) {\,\,$x$};
\node (xp) at (-4,-0.5) {};
\node[above] at (-4,-0.5) {\,\,\,\,$x'$};
\node (x) at (v12) {};
\draw[red,ultra thick] (v11.center) edge (v1.center);
\node (li) at (-2.5,-0.5) {};
\draw[green,thick,densely dashed] (v11.center) edge (li.center);
\node (ri) at (-0.5,-0.5) {};
\draw[green,thick,densely dashed] (v2.center) edge (ri.center);
\node[below] at ($(v3)!.2!(v4)$) {$L_2^-$};
\node[below right] at (-7.5,-0.5) {$L_1^-$};
\node[below] at (-4.5,2)  {$L_0^-$};
\node[below] at (2,3) {$L_2^+$};
\node[below] at (4,-0.5) {$L_1^+$};
\node[below] at (5,-4) {$L_0^+$};
\draw[orange,dotted,thick] (3,-5.5) edge (-9,2);
\draw[orange,dotted,thick] (-9.5,-4.75) edge (4,8);
\path[name path=p2] (3,-5.5) --(-9,2);
\path[name path=34] (v3.center) --(v4.center);
\path[name intersections={of=p2 and 34,by={A}}];

\path[name path=56] (v5) --(v6);
\path[name path=or] (-9.5,-4.75)--(4,8);
\path[name intersections={of=or and 56,by={B}}];

\node (Q) at ($(v5)!.06!(v6)$) {};
\node[above right] at ($(v5)!.06!(v6)$) {$q_0(x)$};
\node (Qp) at ($(v5)!.01!(v6)$) {};
\node[above] at ($(v5)!.01!(v6)$) {$q_0(x')$};
\path[name path=uuu] ($(Q)+(Q)-(x)$)--($(x)+(x)-(Q)+(x)-(Q)$);
\path[name path=uuuu] (Qp.center)--($(xp.center)+(xp.center)-(Qp.center)+(xp.center)-(Qp.center)+(xp.center)-(Qp.center)+(xp.center)-(Qp.center)$);
\path[purple,name intersections={of=34 and uuu,by={D}}];
\node[below]  at (D) {$\,\;\,q_2(x)$};
\draw[purple,name intersections={of=34 and uuuu,by={Dp}}];
\node[below]  at (Dp) {$\;\;\;q_2(x')$};

\draw[green,thick,densely dashed] (Q.center) edge (v6);%
\draw[red,ultra thick] (Q.center) edge (v5);
\draw (Q.center) edge (D.center);
\draw (Qp.center) edge (Dp.center);
\draw[green,thick,densely dashed] (v4.center) edge (Dp.center);%
\draw  (v3)--(D);
\draw [green,thick,densely dashed] (D)--(Dp);
\fill[gray]  (v7) rectangle (v9);
\node [blue,regular polygon, regular polygon sides=5,scale=0.3,fill] at (v11) {};
\node [blue,regular polygon, regular polygon sides=5,scale=0.3,fill] at (-4,-0.5) {};

\node [blue,regular polygon, regular polygon sides=5,scale=0.3,fill] at (D) {};
\node [blue,regular polygon, regular polygon sides=5,scale=0.3,fill] at (Dp) {};
\node [blue,regular polygon, regular polygon sides=5,scale=0.3,fill] at (Q) {};
\node [blue,regular polygon, regular polygon sides=5,scale=0.3,fill] at (Qp) {};
\end{tikzpicture}
\caption{The case where $x$~and~$x'$ are wild points and the lines $K_x$~and~$K_{x'}$ do not intersect.
We argue that the case depicted in the figure cannot occur.
To see this, we see that each point between $q_0(x')$~and~$q_0(x)$ is mapped to a value smaller than~$f(x')$, contradicting the definition of~$q_0(x')$. As we are in the case where $q_0(x')\not=q_0(x)$, we see that $q_0(x')$ is between $q_0(x)$~and~$L_0^+$.}\label{position of lines}
\endgroup
\end{figure}
\end{pf}

The proof above seems to indicate that only the left-hand side ($L_0^- \cup L_1^- \cup L_2^-$) is significant. On the other hand, we have essentially used the fact that the right-hand side of $S$ is unbounded.

We shall now argue that the assertion of the theorem above cannot be improved much and moreover there may be no wild points.
Namely, let $\Err \cdot \Err$ denote the lexicographic square of two copies of $\Err$.
Specifically, $(x,y) < (x',y')$ if and only if either $x < x'$ or else $x = x'$ and $y < y'$.
Then the identity mapping from the Euclidean plane $\Err^2$ onto $\Err \cdot \Err$ is monotone, showing that the plane can actually be mapped to a sufficiently large linearly ordered set.

Now, consider again a set $S \subs \Err^2$ consisting of three pairwise non-parallel lines.
Rotating the plane, if necessary, we may assume that one of these lines is vertical, say, $\sn0 \times \Err$.
Then the identity mapping is monotone and maps $S$ onto the lexicographic sum
\begin{equation*}
  ((-\infty,0) \cdot 2) + (\sn0 \times \Err) + ((0,+\infty) \cdot 2).
\end{equation*}
This in turn embeds into $D = \Err \cdot 2$, the linearly ordered set often called the \emph{double arrow}. Note that $D$ is separable in the interval topology, therefore it has no uncountable family of pairwise disjoint open intervals.
In particular, no open convex planar set admits a one-to-one mapping into $D$.

It should be possible to show that a set consisting of $n$ lines always admits a one-to-one monotone mapping into $\Err \cdot (n-1)$. How about the case where some of the lines are parallel?

The last result shows that admitting a one-to-one monotone map into the real line is not a local property.

\begin{tw}
	Assume $S \subs \Err^2$ consists of three bounded line segments intersecting at a single point.
	Then there exists a one-to-one monotone mapping $\map f S \Err$.
\end{tw}

\begin{pf}[Sketch of proof]
	The construction is visualized in Figure~\ref{FigTrziOdcinkis}. We may assume that one of the segments, say $J$, is (almost) vertical. The infinitely many lines from its two end-points divide the two segments on the left-hand side of $J$ into infinitely many intervals (half-open, half-closed). We map each of them onto an interval in the real line and the enumeration determines the ordering. After that, we map $J$ onto a real interval above all the previous ones. We repeat the same procedure symmetrically on the right-hand side of $J$.
\end{pf}

\begin{figure}
	\centering
\tikzset{
  rgb color/.code args={#1=#2}{%
    \definecolor{tempcolor-#1}{rgb}{#2}%
    \tikzset{#1=tempcolor-#1}%
  },
}
    \iffarben
    \else
\colorlet{L1u}{black}
\colorlet{L2u}{black}
\colorlet{L3u}{black}
\colorlet{L4u}{black}
\colorlet{L5u}{black}
\colorlet{L6u}{black}
\colorlet{L7u}{black}
\colorlet{L8u}{black}
\colorlet{L9u}{black}
\colorlet{L10u}{black}
\colorlet{L11u}{black}
\colorlet{L12u}{black}
\colorlet{L13u}{black}
\colorlet{L14u}{black}
\colorlet{L15u}{black}
\colorlet{L16u}{black}
\colorlet{L17u}{black}

\colorlet{L1d}{black}
\colorlet{L2d}{black}
\colorlet{L3d}{black}
\colorlet{L4d}{black}
\colorlet{L5d}{black}
\colorlet{L6d}{black}
\colorlet{L7d}{black}
\colorlet{L8d}{black}
\colorlet{L9d}{black}
\colorlet{L10d}{black}
\colorlet{L11d}{black}
\colorlet{L12d}{black}
\colorlet{L13d}{black}
\colorlet{L14d}{black}
\colorlet{L15d}{black}
\colorlet{L16d}{black}
\colorlet{L17d}{black}

\colorlet{R1u}{black}
\colorlet{R2u}{black}
\colorlet{R3u}{black}
\colorlet{R4u}{black}
\colorlet{R5u}{black}
\colorlet{R6u}{black}
\colorlet{R7u}{black}
\colorlet{R8u}{black}
\colorlet{R9u}{black}
\colorlet{R10u}{black}
\colorlet{R11u}{black}
\colorlet{R12u}{black}
\colorlet{R13u}{black}
\colorlet{R14u}{black}
\colorlet{R15u}{black}
\colorlet{R16u}{black}
\colorlet{R17u}{black}
\colorlet{R18u}{black}

\colorlet{R1d}{black}
\colorlet{R2d}{black}
\colorlet{R3d}{black}
\colorlet{R4d}{black}
\colorlet{R5d}{black}
\colorlet{R6d}{black}
\colorlet{R7d}{black}
\colorlet{R8d}{black}
\colorlet{R9d}{black}
\colorlet{R10d}{black}
\colorlet{R11d}{black}
\colorlet{R12d}{black}
\colorlet{R13d}{black}
\colorlet{R14d}{black}
\colorlet{R15d}{black}
\colorlet{R16d}{black}
\colorlet{R17d}{black}
\colorlet{R18d}{black}
\colorlet{R19d}{black}

\colorlet{M1u}{black}
\colorlet{M1d}{black}
\fi
\begin{tikzpicture}[scale=.06]
  \useasboundingbox (0, 0) rectangle (240.082, 282.602);
  \draw[black, line width=0.398]
    (6.6538, 162.562)
     -- (148.388, 6.6542);
  \draw[black, line width=0.398]
    (6.6538, 77.5212)
     -- (120.041, 275.949);
  \draw[black, line width=0.398]
    (148.388, 6.6542)
     -- (52.0109, 156.892);
  \draw[black, line width=0.398]
    (148.388, 6.6542)
     -- (75.8332, 153.914);
  \draw[black, line width=0.398]
    (120.041, 275.949)
     -- (58.3608, 105.681);
  \node[font=\fontsize{9.962}{11.9544}\selectfont, text=black]
     at (28.603, 163.026) {1};
  \draw[L1u, line width=1.594]
    (6.6538, 162.562)
     -- (52.0109, 156.892);
  \draw[black, line width=0.398]
    (10.1736, 79.437)
     -- (58.3608, 105.681);
  \node[text=black]
     at (36.5,88.5) {2};
  \draw[black, line width=0.398]
    (52.0109, 156.892)
     -- (75.8332, 153.914);
  \node[text=black]
     at (61.432, 158.922) {3};
  \draw[L2u,line width=1.594]
    (52.0109, 156.892)
     -- (75.8332, 153.914);
  \node[text=black]
     at (68.5,105.5) {4};
  \draw[black, line width=0.398]
    (119.114, 272.429)
     -- (78.0018, 116.378);
  \draw[L3u, line width=1.594]
    (75.8332, 153.914)
     -- (87.5079, 152.455);
  \node[text=black]
     at (80.011, 156.704) {5};
  \draw[line width=1.594,rgb color={draw=0.502 1 0.502}, line width=0.398]
    (78.0018, 116.378)
     -- (90.8737, 123.388);
  \node[text=black]
     at (87.5,114) {6};
  \draw[white, line width=0.0]
    (120.041, 141.302) rectangle (120.041, 141.302);
  \draw
    (90.8732, 123.388)
     -- (120.041, 275.949);
  \draw
    (87.5059, 152.453)
     -- (148.388, 6.6542);
  \draw
    (98.0178, 127.279)
     -- (120.041, 275.949);
  \draw
    (96.2217, 151.363)
     -- (148.388, 6.6542);
  \draw
    (101.488, 150.705)
     -- (148.388, 6.6542);
  \draw
    (103.774, 130.414)
     -- (120.041, 275.949);
  \draw
    (107.443, 132.413)
     -- (120.041, 275.949);
  \draw
    (105.979, 150.143)
     -- (148.388, 6.6542);
  \draw
    (108.967, 149.77)
     -- (148.388, 6.6542);
  \draw
    (110.696, 134.184)
     -- (120.041, 275.949);
  \draw
    (112.925, 135.398)
     -- (120.041, 275.949);
  \draw
    (111.701, 149.428)
     -- (148.388, 6.6542);
  \draw
    (113.624, 149.187)
     -- (148.388, 6.6542);
  \draw[L4u, line width=1.594]
    (87.5065, 152.455)
     -- (96.2217, 151.363);
  \draw[black]
    (120.041, 275.949)
     -- (205.694, 137.676);
  \draw[black]
    (197.145, 181.265)
     -- (148.388, 6.6542);
  \draw[black]
    (148.388, 6.6542)
     -- (183.347, 173.751)
     -- (183.347, 173.751);
  \draw[black]
    (185.672, 140.179)
     -- (120.041, 275.949);
  \draw[black]
    (120.041, 275.949)
     -- (176.562, 141.318);
  \draw[black]
    (148.388, 6.6542)
     -- (166.743, 164.707);
  \draw[black]
    (164.206, 142.863)
     -- (120.041, 275.949);
  \draw[black]
    (148.388, 6.6542)
     -- (158.455, 160.194);
  \draw[black]
    (157.374, 143.717)
     -- (120.041, 275.949)
     -- (120.041, 275.949);
  \draw[black]
    (148.388, 6.6542)
     -- (153.486, 157.488);
  \draw[L5u, line width=1.594]
    (96.2217, 151.363)
     -- (101.488, 150.705);
  \draw[L6u,line width=1.594]
    (101.488, 150.705)
     -- (105.979, 150.143);
  \draw[line width=1.594,rgb color={draw=1 0 0.502}]
    (105.979, 150.143)
     -- (133.648, 146.684);
  \draw[L7u,line width=1.594]
    (108.966, 149.769)
     -- (105.979, 150.143);
  \draw[L8u, line width=1.594]
    (108.966, 149.769)
     -- (111.701, 149.427);
  \draw[L9u, line width=1.594]
    (111.701, 149.427)
     -- (113.623, 149.187);
  \draw[L10u, line width=1.594]
    (113.624, 149.187)
     -- (133.648, 146.683);
  \draw[L1d, line width=1.594]
    (6.6538, 77.5212)
     -- (58.3627, 105.682);
  \draw[L4d,line width=1.594]
    (90.8732, 123.388)
     -- (98.0178, 127.279);
  \draw[L5d, line width=1.594]
    (98.0178, 127.279)
     -- (103.774, 130.415);
  \draw[L6d, line width=1.594]
    (103.774, 130.415)
     -- (107.443, 132.414);
  \draw[L7d, line width=1.594]
    (107.443, 132.414)
     -- (112.925, 135.4);
  \draw[L8d, line width=1.594]
    (110.696, 134.186)
     -- (112.925, 135.4);
  \draw[L9d, line width=1.594]
    (112.925, 135.4)
     -- (133.648, 146.683);
  \draw[black]
    (166.743, 164.707)
     -- (197.145, 181.265);
  \draw[R9u,line width=1.594]
    (183.347, 173.75)
     -- (197.145, 181.265);
  \draw[R9d,line width=1.594]
    (205.694, 137.676)
     -- (185.672, 140.179);
  \draw[R8d,line width=1.594]
    (176.562, 141.318)
     -- (185.672, 140.179);
  \draw[R8u,line width=1.594]
    (183.347, 173.75)
     -- (172.341, 167.756);
  \draw[black]
    (172.341, 167.756)
     -- (148.388, 6.6542);
  \draw[black]
    (168.559, 142.318)
     -- (120.041, 275.949);
  \draw[black]
    (161.473, 161.837)
     -- (148.388, 6.6542);
  \draw[black]
    (159.918, 143.399)
     -- (120.041, 275.949);
  \draw[black]
    (155.371, 158.514)
     -- (148.388, 6.6542);
  \draw[black]
    (154.706, 144.05)
     -- (120.041, 275.949);
  \draw[R2u,line width=1.594]
    (151.464, 156.386)
     -- (153.486, 157.487);
  \draw[R1u,line width=1.594]
    (133.648, 146.683)
     -- (151.464, 156.386);
  \draw[R6u,line width=1.594]
    (161.473, 161.837)
     -- (166.743, 164.707);
  \draw[R5u,line width=1.594]
    (158.455, 160.193)
     -- (161.473, 161.837);
  \draw[R4u,line width=1.594]
    (155.371, 158.514)
     -- (158.455, 160.193);
  \draw[R3u,line width=1.594]
    (153.486, 157.487)
     -- (155.371, 158.514);
  \draw[R7d,line width=1.594]
    (168.559, 142.318)
     -- (176.562, 141.319);
  \draw[R6d,line width=1.594]
    (164.206, 142.863)
     -- (168.559, 142.318);
  \draw[R5d,line width=1.594]
    (159.918, 143.399)
     -- (164.206, 142.863);
  \draw[R4d,line width=1.594]
    (157.375, 143.717)
     -- (159.918, 143.399);
  \draw[R3d,line width=1.594]
    (157.375, 143.717)
     -- (154.706, 144.05);
  \draw[R2d,line width=1.594]
    (154.706, 144.05)
     -- (153.039, 144.259);
  \draw[black]
    (153.039, 144.259)
     -- (133.648, 146.683);
  \draw[R7u,line width=1.594]
    (166.743, 164.707)
     -- (172.341, 167.756);
  \draw[L2d, line width=1.594]
    (58.3616, 105.6833)
     -- (78.0008, 116.3775);
  \draw[L3d, line width=1.594]
    (78.0018, 116.378)
     -- (90.8734, 123.3878);
  \draw[R1d,line width=1.594]
    (133.648, 146.683)
     -- (153.0389, 144.259);
  \draw[M1u,line width=1.594]
    (120.041, 275.949)
     -- (133.648, 146.683);
  \draw[M1d,line width=1.594]
    (133.648, 146.683)
     -- (148.388, 6.6542);
\end{tikzpicture}
	\caption{Three segments intersecting in a point can be mapped with a one-to-one monotone mapping into~$\Err$.
The figure takes care about the definition of the mapping on the left-hand side. We partition the left sides of the segments into countably many half-open intervals as depicted and order their images as indicated by the numbers.}\label{FigTrziOdcinkis}
\end{figure}

\section{Final remarks}

Our main result, namely, a classification of monotone mappings (Theorem~\ref{THMgdyifgi}, holds when the real plane is replaced by $K^2$, where $K$ is a subfield of $\Err$. Indeed, the ``regular'' case relies on Theorem~\ref{THMhomographies}, whose existence part is purely algebraic (the homography would have coefficients in $K$) and the uniqueness arguments are exactly the same, as $K^2$ is dense in $\Err^2$ and moreover every line segment in $K^2$ is dense in its closure in $\Err^2$. The ``pathological'' cases are proved in exactly the same way.

It is not clear what happens if $\Err$ is replaced by an arbitrary ordered field. Homographies are as before, however, the density arguments might be more subtle, if the field is non-archimedean.

Another possible direction of future research is a classification of monotone mappings from convex sets in higher dimensions. There might be more ``pathological'' cases here, and it is not clear to us how to classify them.

Yet another aspect, worth studying in our opinion, is which ``simple'' planar sets admit a one-to-one monotone mapping into the real line. This problem formally depends on the betweenness relation only, however, we have seen that it is not a local property and some small perturbations of the set (e.g.\ parallel vs.\ non-parallel lines) may give different results.

\section{Appendix: A proof of Proposition~\ref{propreducn2}}

First of all, recall that every separable infinite-dimensional space has linear dimension over $\Err$ precisely the continuum (see~\cite{Lac73} for a simple proof).
Thus, we may assume that $X \subs V$, where $V$ is, say, the separable infinite-dimensional Hilbert space. In this setting, we shall show that there exists a bounded linear operator satisfying the assertion of Proposition~\ref{propreducn2}.

For each triple $s = \{a,b,c\} \subs X$ consisting of non-collinear points, we define $W_s$ to be the set of all linear maps onto $V_0$ that make the points $a,b,c$ collinear. This is a closed linear subspace of the space $L(V,V_0)$ of all bounded linear maps from $V$ onto $V_0$. Its co-dimension is $2$.
Let $W$ denote the union of all $W_s$, where $s$ is as above.
It is necessary and sufficient to show that $W \ne L(V,V_0)$.
This is evidently true if $X$ is countable, due to the Baire category theorem.

In general, one can use an argument from logic. Namely, the linear space $W$ has a concrete definition with parameters in $X$ and in the model of set theory obtained by collapsing the continuum to $\aleph_1$ the set $L(V,V_0) \setminus W$ is nonempty. Hence the same set is nonempty in the original model we are working with. This argument uses Gödel's completeness, namely, a sentence is a theorem if and only if it holds in every model of set theory (see~\cite{Hermes}, Chapter~V).

Note that every Banach space of density $\aleph_1$ can be covered by $\aleph_1$ many closed hyperplanes. This is because such a space is the union of an $\aleph_1$-chain of closed separable subspaces and each of them can be enlarged (by the Hahn--Banach Theorem) to a closed hyperplane. This indicates that perhaps the linear map required by Proposition~\ref{propreducn2} might be unbounded, if we take $V$ to be a non-separable normed space. On the other hand, the statement has nothing to do with topology, even though the proof above gives a bounded linear map.

\paragraph{Acknowledgements.} The research of the first author was supported by grant 20-22230L (Czech Science Foundation) and of the second and third author by the University of Silesia Mathematics Department (Iterative Functional Equations and Real Analysis program).
\FloatBarrier
\bibliographystyle{alpha}
\bibliography{MonotoneMappings}
\end{document}